\documentclass{amsart}
\usepackage{graphicx}
\usepackage{xypic}
\usepackage{psfrag}

\theoremstyle{plain}
\newtheorem{prop}{Proposition}
\newtheorem{thm}{Theorem}
\newtheorem{example}{Example}

\newtheorem*{cor}{Corollary}
\newtheorem{lemma}{Lemma}

\theoremstyle{remark}
\theoremstyle{definition}
\newtheorem*{remark}{Remark}

\newcommand{\refthm}[1]{Theorem~\ref{#1}}
\newcommand{\refprop}[1]{Proposition~\ref{#1}}
\newcommand{\reflemma}[1]{Lemma~\ref{#1}}

\newcommand{\Schub}{{\mathfrak S}}
\newcommand{\Groth}{{\mathfrak G}}

\newcommand{\Z}{{\mathbb Z}}

\newcommand{\C}{{\mathbb C}}
\newcommand{\cO}{{\mathcal O}}

\DeclareMathOperator{\Hom}{Hom}

\DeclareMathOperator{\GL}{GL}

\newcommand{\pic}[2]{\includegraphics[scale=0.#1]{#2}}

\newcommand{\ignore}[1]{}

\newcommand{\bu}{{\mathbf u}}
\newcommand{\bw}{{\mathbf w}}

\DeclareMathOperator{\tp}{Tp}
\newcommand{\wb}{\overline}

\psfrag{_B1}{$b_1$}
\psfrag{_B2}{$b_2$} 
\psfrag{_B3}{$b_3$}
\psfrag{Ri1i}{$R_{i-1,i}$}
\psfrag{Ri1i1}{$R_{i-1,i+1}$}
\psfrag{Ri1n}{$R_{i-1,n}$}
\psfrag{dots}{$\cdots$}

\begin{document}

\title[Positivity of quiver coefficients]{Positivity of quiver
  coefficients \\ through Thom polynomials}

\author{Anders~S.~Buch}
\address{Matematisk Institut, Aarhus Universitet, Ny Munkegade, 8000
  {\AA}rhus C, Denmark}
\email{abuch@imf.au.dk}
\author{L\'aszl\'o~M.~Feh\'er}
\address{Department of Analysis, Eotvos University, Budapest, Hungary}
\email{lfeher@math-inst.hu}
\author{Rich\'ard Rim\'anyi}
\address{Department of Mathematics, The University of North Carolina
  at Chapel Hill, USA}
\email{rimanyi@email.unc.edu}

\date{November 11, 2003}
\thanks{Supported by FKFP0055/2001 (2nd and 3rd author)}

\maketitle


\section{Introduction}

Let $(e_0,e_1,\dots,e_n)$ be a dimension vector of non-negative
integers.  The space $V = \Hom(\C^{e_0},\C^{e_1}) \oplus \cdots \oplus
\Hom(\C^{e_{n-1}},\C^{e_n})$ of equioriented quiver representations of
type $A$ has a natural action of the group $G = \GL(e_0) \times \dots
\times \GL(e_n)$ given by $(g_0,\dots,g_n) . (\phi_1,\dots,\phi_n) =
(g_1 \phi_1 g_0^{-1}, \dots, g_n \phi_n g_{n-1}^{-1})$.  An orbit
$r$ of this action is characterized by its set of {\em rank
  conditions\/} $\{ r_{ij} \}$ for $0 \leq i < j \leq n$, where
$r_{ij}$ is the rank of the composed map $\phi_j \phi_{j-1} \cdots
\phi_{i+1}$ for any point in this orbit.  In this paper we study the
$G$-equivariant cohomology class of the orbit closure $\overline r$.
We will call this class for the {\em Thom polynomial\/} of the orbit,
and we denote it by $\tp_r$.

This Thom polynomial can be regarded as a formula for the degeneracy
locus obtained when the integers $r_{ij}$ are used as rank conditions
for a sequence of vector bundles and bundle maps (see e.g.\ 
\cite{fulton:notes} for the translation).  Buch and Fulton gave a
formula expressing the cohomology class of such a degeneracy locus as
a linear combination of products of Schur determinants
\cite{buch.fulton:chern}.  When interpreted for Thom polynomials, this
formula has the form
\begin{equation} \label{eqn:orig_qf}
  \tp_r = \sum_\lambda c_\lambda(r) \, s_{\lambda_1}(x^1;x^0) \,
   s_{\lambda_2}(x^2;x^1) \cdots s_{\lambda_n}(x^n;x^{n-1})
\end{equation}
where the sum is over certain sequences of partitions $\lambda$, and
the symbol $x^i$ denotes the Chern roots $\{ x^i_1, \dots, x^i_{e_i}
\}$ of the $i$'th factor of $G$.  The {\em quiver coefficients\/}
$c_\lambda(r)$ appearing in this formula are integers uniquely
determined by (\ref{eqn:orig_qf}) in addition to the condition that
$c_\lambda(r) = c_\lambda(r+k)$ for all $k \geq 0$, where $r+k$
denotes the rank conditions $\{ r_{ij} + k \}$ obtained by adding the
integer $k$ to the original rank conditions.  Although the formula for
quiver coefficients in \cite{buch.fulton:chern} does not reveal their
signs, it was conjectured that all quiver coefficients are
non-negative.

Feh\'er and Rim\'anyi suggested a different method for computing Thom
polynomials in \cite{feher.rimanyi:calculation, feher.rimanyi:thom},
which works more generally for all quiver representations associated
to Dynkin diagrams.  In this approach, the Thom polynomial $\tp_r$ is
obtained as the unique solution to a system of linear equations.

The $G$-orbits in the representation space $V$ were first classified
by Abeasis and Del Fra using {\em lace diagrams\/}
\cite{abeasis.del-fra:degenerations}.  An important idea in recent
work of Knutson, Miller, and Shimozono \cite{knutson.miller.ea:four}
was to reinterpret these lace diagrams as sequences of permutations,
which can be identified with the components of a {\em Gr\"obner
  degeneration\/} of the orbit closure.  In a talk about this work
given by E.~Miller at the Boston AMS-meeting in October 2002, the
following {\em component formula\/} was conjectured, which expresses
the Thom polynomial $\tp_r$ as a sum of products of Schubert
polynomials:
\[ \tp_r = \sum_{(w_1,\dots,w_n)} \Schub_{w_1}(x^1;x^0) \,
   \Schub_{w_1}(x^2;x^1) \cdots \Schub_{w_n}(x^n;x^{n-1}) \,.
\]
This sum is over all {\em minimal\/} lace diagrams, whose definition
is recalled in section \ref{sec:component}.  This conjecture was
subsequently proved independently by the authors of
\cite{knutson.miller.ea:four} and the third author of the present
paper.  The main goal of this paper is to present the Hungarian
approach, which consists of simply verifying that the component
formula satisfies the required equations for being a Thom polynomial.

The component formula also has a stable variant, where the Schubert
polynomials are replaced with {\em Stanley symmetric functions}.  This
version of the formula was first proved in
\cite{knutson.miller.ea:four}.  Since Stanley symmetric functions are
Schur positive \cite{edelman.greene:balanced,
  lascoux.schutzenberger:structure}, the stable component formula
implies that quiver coefficients are non-negative.  In this paper we
give a simple argument that the two versions of the component formula
are equivalent, thus obtaining a short proof of the non-negativity of
quiver coefficients based on Thom polynomial theory.  In comparison,
the approach of \cite{knutson.miller.ea:four} relies on two different
geometric constructions, one of which is the above mentioned Gr\"obner
degeneration, and the other being a {\em ratio formula\/} derived from
a geometric study of the Zelevinsky map \cite{zelevinsky:two,
  lakshmibai.magyar:degeneracy}.

Part of our verification of the component formula consists of proving
that this formula is symmetric in each set of variables $x^i$.  This
argument can also be turned around to show that a linear combination
of products of Schubert polynomials over minimal lace diagrams is
symmetric if and only if all coefficients are equal.  This in turn
makes it possible to prove the component formula directly from the
Gr\"obner degeneration, at least up to a constant, which can then be
determined by applying the original quiver formula
\cite{buch.fulton:chern}.  We will explain this alternative proof in
section \ref{sec:grobner}.

We remark that the component formula can also be derived
combinatorially \cite{yong:embedding, buch:alternating} from the ratio
formula of \cite{knutson.miller.ea:four}.  In fact, among the four
geometric approaches to quiver formulas currently known to us
\cite{buch.fulton:chern,feher.rimanyi:thom,knutson.miller.ea:four},
only the original approach of \cite{buch.fulton:chern} (which is based
on resolution of singularities for quiver varieties) offers no easy
path to positivity of quiver coefficients.  On the other hand, the
original approach arguably makes the question of positivity more
natural to ask.

The component formula also has a $K$-theory variant
\cite{buch:alternating, miller:alternating}, which implies that the
$K$-theoretic quiver coefficients defined in \cite{buch:grothendieck}
have alternating signs.  This formula expresses the structure sheaf of
a quiver variety as an alternating sum of products of Grothendieck
polynomials indexed by {\em KMS-factorizations}, which generalize
minimal lace diagrams.  In the last section we apply the methods of
this paper to give a new description of KMS-factorizations based on
transformations of lace diagrams.

This paper is organized as follows.  In section~\ref{sec:component} we
explain basic notions like minimal lace diagrams and Schubert
polynomials, and we prove that the component formula is symmetric and
equivalent to the stable component formula.  In section~\ref{sec:thom}
we prove the component formula using Thom polynomial theory, while
section~\ref{sec:grobner} contains the alternative proof based on the
Gr\"obner degeneration of \cite{knutson.miller.ea:four}.
Section~\ref{sec:ktheory} finally contains the classification of
KMS-factorizations.

We thank P.~Pragacz, A.~Weber, and the Banach Center in Warsaw for
their hospitality while part of this work was carried out.


\section{The component formula}
\label{sec:component}

A {\em lace diagram\/} for the dimension vector $(e_0,\dots,e_n)$ is a
diagram of dots arranged in columns, with $e_i$ dots in column $i$,
together with line segments connecting dots of consecutive columns.
Each dot may be connected to at most one dot in the column to the left
of it, and to at most one dot in the column to the right of it.  The
corresponding orbit $r$ satisfies that the rank condition $r_{ij}$ is
the number of connections from column $i$ to column $j$
\cite{abeasis.del-fra:degenerations}.  For example, the following lace
diagram corresponds to an orbit $r$ of quiver representations through
5 vector spaces of dimensions $(e_0,\dots,e_4) = (3,4,3,3,2)$, and we
have $r_{01} = r_{02} = 2$, $r_{03} = 1$, etc.
\[ \pic{50}{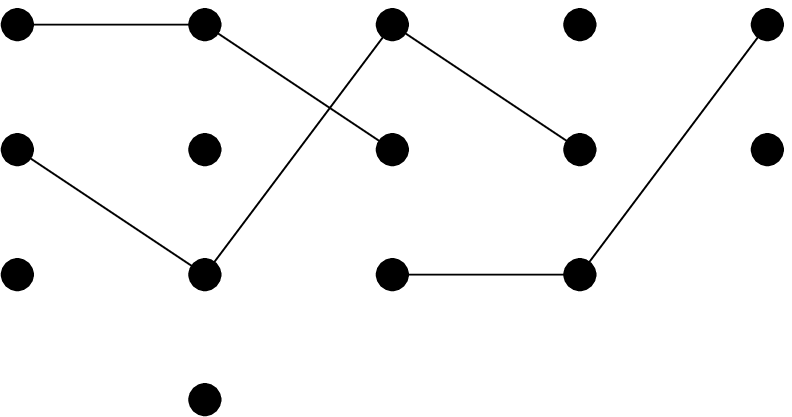} \]

A lace diagram may be identified with a sequence $(w_1,\dots,w_n)$ of
permutations \cite{knutson.miller.ea:four} (see also
\cite{fulton:flags}).  Here we let $w_i$ be the permutation of minimal
length such that $w_i(q) = p$ whenever dot $q$ of column $i$ is
connected to dot $p$ of column $i-1$.  Equivalently, this permutation
describes the connections from the $i$'th to the $i-1$'st column of an
{\em extension\/} of the lace diagram.  This extended diagram is
constructed by adding extra dots to the columns, so that the original
dots without connections to both sides can be connected to the new
dots.  For example, the above lace diagram is extended as follows; in
particular we have $w_2 = 31524$.
\[ \pic{50}{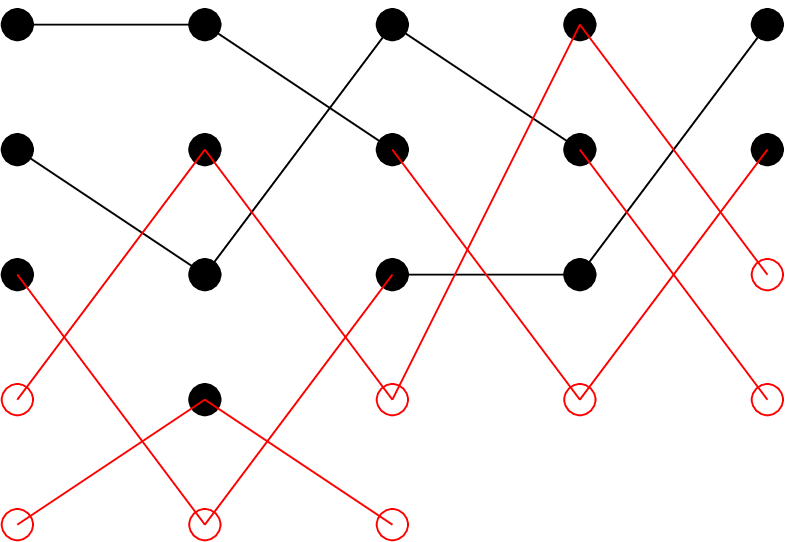} \]

Notice that a sequence $(w_1,\dots,w_n)$ of permutations represent a
lace diagram for the dimension vector $(e_0,\dots,e_n)$ if and only if
each permutation $w_i$ is a {\em partial permutation\/} from $e_i$
elements to $e_{i-1}$ elements, which means that all descent positions
of $w_i$ are smaller than or equal to $e_i$, and all descent positions
of $w_i^{-1}$ are smaller than or equal to $e_{i-1}$.

A {\em strand\/} of a lace diagram is a maximal sequence of connected
dots and line segments, and the {\em extension\/} of a strand is
obtained by also including the extra line segments that it is directly
connected to in the extended lace diagram.

The {\em length\/} of a lace diagram is the sum of the lengths of the
permutations $w_i$, or equivalently the total number of crossings in
the extended lace diagram.  Notice that the smallest possible length
of a lace diagram for an orbit $r$ is equal to the codimension $d(r) =
\sum_{i < j} (r_{i,j-1}-r_{ij})(r_{i+1,j}-r_{ij})$ of the orbit.  This
follows because all of the $r_{i+1,j}-r_{ij}$ strands starting in
column $i+1$ and passing through column $j$ must intersect all of the
$r_{i,j-1}-r_{ij}$ strands passing through column $i$ and terminating
in column $j-1$.  The lace diagram is called {\em minimal\/} if its
length is equal to $d(r)$.  This is equivalent to demanding that (the
extensions of) any two strands can cross at most once, and not at all if
they start or end at the same column
\cite[Thm.~3.8]{knutson.miller.ea:four}.

To state the component formula, we also need the Schubert polynomials
of Lascoux and Sch\"utzenberger
\cite{lascoux.schutzenberger:structure}.  The {\em divided difference
  operator\/} $\partial_{a,b}$ with respect to two variables $a$ and
$b$ is defined by
\[ \partial_{a,b}(f) = \frac{f(a,b) - f(b,a)}{a-b} \,, \]
where $f$ is any polynomial in these (and possibly other) variables.
The {\em double Schubert polynomials\/} $\Schub_w(x;y) =
\Schub_w(x_1,\dots,x_m;y_1,\dots,y_m)$ given by permutations $w \in
S_m$ are uniquely determined by the identity
\begin{equation} \label{eqn:schub_dd}
  \partial_{x_i,x_{i+1}}(\Schub_w(x;y)) = \begin{cases}
  \Schub_{w s_i}(x;y) & \text{if $w(i) > w(i+1)$} \\
  0 & \text{if $w(i) < w(i+1)$}
  \end{cases} 
\end{equation}
together with the expression 
\[ \Schub_{w_0}(x;y) = \prod_{i+j \leq m} (x_i - y_j) \]
for the longest permutation $w_0$ in $S_m$.  Using that $\Schub_w(y;x)
= (-1)^{\ell(w)} \Schub_{w^{-1}}(x;y)$ we similarly have that
$\partial_{y_i,y_{i+1}}(\Schub_w(x;y))$ is equal to $-\Schub_{s_i
  w}(x;y)$ if $\ell(s_i w) < \ell(w)$, and is zero otherwise.  If $k$
and $l$ are the last descent positions of $w$ and $w^{-1}$,
respectively, then only the variables $x_1,\dots,x_k,y_1,\dots,y_l$
occur in $\Schub_w(x;y)$.

The component formula can now be stated as follows.  Recall that the
variables $x^i_j$ are the Chern roots of the group $G$ of the
introduction.

\begin{thm} \label{thm:ns_comp}
The Thom polynomial of a $G$-orbit $r$ is given by
\[ \tp_r = \sum_{(w_1,\dots,w_n)} \Schub_{w_1}(x^1;x^0) \,
   \Schub_{w_2}(x^2;x^1) \cdots \Schub_{w_n}(x^n;x^{n-1})
\]
where the sum is over all minimal lace diagrams representing the
orbit.
\end{thm}

It follows from this theorem that the component formula is symmetric
in each set of variables $x^i$.  This can also be proved directly.  We
let
\begin{equation} \label{eqn:comppoly}
  Q_r = \sum_{(w_1,\dots,w_n)} \Schub_{w_1}(x^1;x^0) \,
  \Schub_{w_2}(x^2;x^1) \cdots \Schub_{w_n}(x^n;x^{n-1})
\end{equation}
denote the polynomial of the component formula.

\begin{lemma} \label{lemma:comp_sym}
  The polynomial $Q_r$ is symmetric in each set of variables $x^i$.
\end{lemma}
\begin{proof}
  We must show that for any $0 \leq i \leq n$ and $1 \leq j < e_i$,
  the divided difference operator $\partial^i_j =
  \partial_{x^i_j,x^i_{j+1}}$ maps $Q_r$ to zero.  Notice at first
  that any minimal lace diagram $(w_1,\dots,w_n)$ must satisfy that
  $\ell(s_k w_1) > \ell(w_1)$ for $k < e_0$ and $\ell(w_n s_k) >
  \ell(w_n)$ for $k < e_n$.  Using (\ref{eqn:schub_dd}) this implies
  that $\partial^i_j(Q_r) = 0$ for $i=0$ or $i=n$.
  
  Given any sequence of permutations $(w_1,\dots,w_n)$ we write
  $\Schub(w_1,\dots,w_n) = \prod \Schub_{w_i}(x^i;x^{i-1})$ for the
  corresponding product of Schubert polynomials.  Now suppose that $1
  \leq i \leq n-1$ and let $(w_1,\dots,w_n)$ be a minimal lace diagram
  for $r$.  There are four cases to consider:
  
  (i) $w_i(j) < w_i(j+1)$ and $w_{i+1}^{-1}(j) < w_{i+1}^{-1}(j+1)$.  We
  get $\partial^i_j(\Schub(w_1,\dots,w_n)) = 0$.

  (ii) $w_i(j) < w_i(j+1)$ and $w_{i+1}^{-1}(j) > w_{i+1}^{-1}(j+1)$.
  We get $\partial^i_j(\Schub(w_1,\dots,w_n)) =
  - \Schub(w_1,\dots,w_i,s_j w_{i+1},\dots,w_n)$.
  
  (iii) $w_i(j) > w_i(j+1)$ and $w_{i+1}^{-1}(j) < w_{i+1}^{-1}(j+1)$.
  We get $\partial^i_j(\Schub(w_1,\dots,w_n)) = \Schub(w_1,\dots,w_i
  s_j, w_{i+1},\dots,w_n)$.

  (iv) $w_i(j) > w_i(j+1)$ and $w_{i+1}^{-1}(j) > w_{i+1}^{-1}(j+1)$.
  This is impossible since $(w_1,\dots,w_i s_j, s_j w_{i+1}, \dots,
  w_n)$ would be a shorter lace diagram for the orbit $r$.
  
  Notice that if our minimal lace diagram $\bw = (w_1,\dots,w_n)$
  falls in one of the cases (ii) or (iii), then the sequence $\bw' =
  (w_1,\dots,w_i s_j, s_j w_{i+1}, \dots, w_n)$ is also a minimal lace
  diagram for $r$.  For example, if $w_i(j) > w_i(j+1)$ then since two
  crossing strands cannot both terminate at column $i$, we must have
  $w_{i+1}^{-1}(j) \leq e_{i+1}$, which implies that $\bw'$ is also a
  lace diagram.  Since $\partial^i_j(\Schub(\bw) + \Schub(\bw')) = 0$,
  we conclude that $\partial^i_j(Q_r) = 0$ as required.
\end{proof}

The {\em double Stanley symmetric function\/} $F_w$ for a permutation
$w$ is defined by
\[ F_w(x_1,\dots,x_p;y_1,\dots,y_q) = 
   \Schub_{1^k \times w}(x_1,\dots,x_p,0,\dots,0; 
   y_1,\dots,y_q,0,\dots,0) \,,
\]
where $k$ is any integer larger than $p$ and $q$, and the shifted
permutation $1^k \times w$ acts as the identity on the set
$\{1,\dots,k\}$, and maps $k+j$ to $k+w(j)$ for $j \geq 1$.  We also
need the identity
\begin{equation} \label{eqn:schubspec}
  \Schub_{1^k \times w}(0^k,x_1,\dots,x_m ; 0^k,y_1,\dots,y_m) =
  \Schub_w(x_1,\dots,x_m; y_1,\dots,y_m)
\end{equation}
where $0^k$ denotes $k$ zeros.  This identity is proved in
\cite[Cor.~4]{buch.rimanyi:specializations}.

The following consequence of \reflemma{lemma:comp_sym} shows that
\refthm{thm:ns_comp} is equivalent to the stable component formula,
which states that the Thom polynomial $\tp_r$ equals the sum of
products of Stanley symmetric functions in the corollary.  By the
Schur positivity of Stanley symmetric functions
\cite{edelman.greene:balanced, lascoux.schutzenberger:structure}, this
formula implies that quiver coefficients are non-negative.  The
statement of the corollary was first proved in
\cite{knutson.miller.ea:four} using a combination of geometry and
combinatorics.

\begin{cor}[Knutson, Miller, Shimozono]
  For any orbit $r$ we have
\[ Q_r = \sum_{(w_1,\dots,w_n)} F_{w_1}(x^1;x^0) \, F_{w_2}(x^2;x^1)
   \cdots F_{w_n}(x^n;x^{n-1})
\]
where the sum is over all minimal lace diagrams for $r$.
\end{cor}
\begin{proof}
  Let $r+k$ be denote the orbit corresponding to the dimension vector
  $(e_0+k,\dots,e_n+k)$ and rank conditions $\{ r_{ij}+k \}$.  The
  above discussion of lace diagrams implies that the minimal lace
  diagrams for $r+k$ are exactly those obtained by adding $k$ strands
  of length $n$ to the top of a minimal lace diagram for $r$
  \cite[Cor.~4.12]{knutson.miller.ea:four}.  Equivalently, such a
  diagram is given by a sequence of permutations $(1^k \times w_1,
  \dots, 1^k \times w_n)$, for which $(w_1,\dots,w_n)$ is a minimal
  lace diagram for $r$.  For $k \geq \max(e_0,\dots,e_n)$, the
  symmetry of the polynomial $Q_{r+k}$ therefore implies that
\begin{multline*} 
 \sum_{\bw} \prod_{i=1}^n \Schub_{w_i}(x^i;x^{i-1})
 = \sum_{\bw} \prod_{i=1}^n 
   \Schub_{1^k \times w_i}(0^k,x^i \, ;\,  0^k,x^{i-1}) \\
 = \sum_{\bw} \prod_{i=1}^n 
   \Schub_{1^k \times w_i}(x^i,0^k \, ;\,  x^{i-1},0^k)
 = \sum_{\bw} \prod_{i=1}^n F_{w_i}(x^i;x^{i-1})
\end{multline*}
where the sums are over all minimal lace diagrams $\bw =
(w_1,\dots,w_n)$ for $r$.
\end{proof}

We remark that the first equality in the above proof can also be
deduced from \refthm{thm:ns_comp} together with the property
$c_\lambda(r) = c_\lambda(r+k)$ of quiver coefficients.  Since the
proof of \refthm{thm:ns_comp} using Thom polynomial theory in section
\ref{sec:thom} does not rely on the corollary, one can therefore prove
that quiver coefficients are non-negative without relying on
(\ref{eqn:schubspec}).  However, the alternative proof of the
component formula in section \ref{sec:grobner} does rely on the
corollary, which makes the given combinatorial proof preferable.



\section{Proof using Thom polynomials}
\label{sec:thom}

Let $G$ be a complex Lie group acting on a vector space $V$ with
finitely many orbits.  An orbit $\mu$ of complex codimension $d$ has
an associated $G$-characteristic class $\tp_\mu \in H^{2d}(BG) =
H^{2d}(BG;\Z)$ called its {\em Thom polynomial}.  When $\mu$ is an
orbit of a space of quiver representations as in the introduction,
this Thom polynomial is equivalent to the quiver formula
(\ref{eqn:orig_qf}).

We let $G_\mu$ denote the stabilizer subgroup of a point
$p_\mu$ in $\mu$.  The inclusion of $G_\mu$ into $G$ induces a map
$BG_\mu \to BG$ between the classifying spaces, which gives a ring
homomorphism $\phi_\mu : H^*(BG) \to H^*(BG_\mu)$ on cohomology.  One
can choose a normal slice $N_\mu$ to $\mu$ at $p_\mu$ which is
invariant under the action of the maximal compact subgroup of $G_\mu$.
The Euler class of this action on $N_\mu$ is denoted by $e(\mu) \in
H^*(BG_\mu)$ (note that $H^*(BG_\mu)$ does not change if we pass to
the maximal compact subgroup).

In \cite{feher.rimanyi:calculation} a general theory for computing
Thom polynomials is developed.  The special case of this theory that
is needed here is summarized in the following theorem.

\begin{thm} \label{thm:thom} Let $\mu$ and $\eta$ be orbits of a
  $G$-representation with finitely many orbits.

{\rm (i)} If $\mu\not\subset \overline{\eta}$ then
  $\phi_\mu(Tp_\eta)=0$;

{\rm (ii)} $\phi_\eta(Tp_\eta)=e(\eta)$.

\noindent
Furthermore, if for every orbit $\mu$ the Euler class $e(\mu)$ is not
a zero-divisor in $H^*(BG_\mu)$, then $\tp_\eta$ is uniquely
determined by these conditions.
\end{thm}

For the application to quiver formulas that concerns us here, we use
the group $G = \prod_{i=0}^n \GL(e_i)$ with its usual quiver action on
$V = \bigoplus_{i=1}^n \Hom(\C^{e_{i-1}}, \C^{e_i})$.  In this case the
cohomology ring $H^*(BG)$ is the ring of polynomials in the Chern
roots $x^i_j$, which are symmetric in each group of variables $x^i =
\{x^i_1, \dots, x^i_{e_i} \}$:
\[ H^*(BG) = \Z[ x^i_j \mid 0 \leq i \leq n, 1 \leq j \leq e_i]^{\prod
  S_{e_i}}
\]

In \cite{feher.rimanyi:thom} a combinatorial description of the
cohomology ring $H^*(B G_\mu)$, the restriction map $\phi_\mu$, and
the Euler class $e(\mu)$ was given, which works for representations of
any quiver that is shaped like a Dynkin diagram.  In our case of
equioriented quivers of type A, this works as follows (see \cite[\S
4--5]{feher.rimanyi:thom}).

Let $r \subset V$ be an orbit, and fix a lace diagram $\bw$
representing $r$.  Choose variables $b_1, \dots, b_k$ corresponding to
the strands of $\bw$.  Then $H^*(B G_r) = \Z[b_1,\dots,b_k]$ can be
identified with a polynomial ring in these variables, and $\phi_r :
H^*(B G) \to H^*(B G_r)$ maps each variable $x^i_j$ to the variable of
the strand passing through dot $j$ of column $i$ in $\bw$.  We notice
that this description makes it possible to extend $\phi_r$ to a map on
all polynomials in the Chern roots $x^i_j$.  This extended map depends
on the chosen lace diagram, and is denoted by $\phi_\bw$.  Finally, if
$\bw$ is a minimal lace diagram, then the Euler class $e(r) \in H^*(B
G_r)$ is the product of all differences $(b_p - b_q)$ of variables for
which the extensions of the corresponding strands cross in $\bw$; here
the strand of $b_p$ should have the highest slope at the crossing
point.

\begin{example}
  In the following minimal lace diagram, the strands have been
  labeled with the associated variables.  
\[ \pic{60}{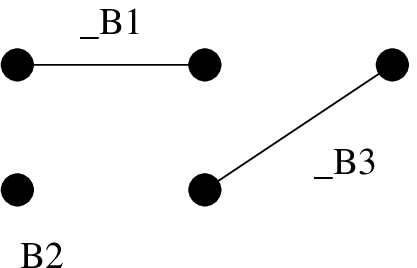} \]
If $r$ denotes the corresponding orbit, then $H^*(BG_r) =
\Z[b_1,b_2,b_3]$, and the map $\phi_r : H^*(BG) \to H^*(BG_r)$ is
given by $\phi_r(x^0_1) = \phi_r(x^1_1) = b_1$, $\phi_r(x^0_2) = b_2$,
and $\phi_r(x^1_2) = \phi_r(x^2_1) = b_3$.  Finally we have $e(r) =
(b_3 - b_1)(b_3-b_2) \in H^*(BG_r)$.
\end{example}


Let $u,w \in S_m$ be permutations.  Our proof of the component formula
uses that the specialization $\Schub_w(b_u;b) =
\Schub_w(b_{u(1)},\dots,b_{u(m)}; b_1,\dots,b_m)$ is zero unless $w
\leq u$ in the Bruhat order on $S_m$, and for $u = w$ we have
\begin{equation} \label{eqn:schub_spec}
  \Schub_u(b_{u(1)},\dots,b_{u(m)}; b_1,\dots,b_m) =
  \prod_{i < j \, ; \, u(i) > u(j)} (b_{u(i)} - b_{u(j)}) \,.
\end{equation}
These statements follow by descending induction on $\ell(w)$ from the
identity
\[ (b_{u(i+1)}-b_{u(i)}) \, \Schub_w(b_u;b) = 
   \Schub_{w s_i}(b_{u s_i}; b) - \Schub_{w s_i}(b_u; b)
\]
which holds whenever $w(i) < w(i+1)$.  The vanishing statement is part
of Goldin's characterization of the Bruhat order
\cite{goldin:cohomology}, and both statements can also be deduced from
\refthm{thm:thom} applied to the representation studied in \cite[\S
4]{ss}.  More general formulas for specializations of Schubert
polynomials are proved in \cite{buch.rimanyi:specializations}.

\begin{proof}[Proof of \refthm{thm:ns_comp} using Thom polynomial
  theory]
  
  \reflemma{lemma:comp_sym} shows that the polynomial $Q_r$ of
  (\ref{eqn:comppoly}) is an element of $H^*(BG)$.  We must show that
  $Q_r$ satisfies the requirements (i) and (ii) of \refthm{thm:thom}.
  
  As in the proof of \reflemma{lemma:comp_sym} we set
  $\Schub(w_1,\dots,w_n) = \prod \Schub_{w_i}(x^i;x^{i-1})$.  Notice
  that if $\bu = (u_1,\dots,u_n)$ is any lace diagram, then
  $\phi_\bu(\Schub(w_1,\dots,w_n))$ is zero unless $w_i \leq u_i$ in
  the Bruhat order for all $i$.  In fact, if $b_1,\dots,b_m$ are the
  variables of the strands through column $i-1$ in the extended lace
  diagram for $\bu$, ordered from top to bottom, then $\phi_\bu$ maps
  the $i$'th factor of $\Schub(w_1,\dots,w_n)$ to
  $\Schub_{w_i}(b_{u_i};b)$.
  
  Now suppose that $s \subset V$ is an orbit which is not contained in
  the closure of $r$.  This implies that $s_{ij} > r_{ij}$ for some $0
  \leq i < j \leq n$.  Choose a lace diagram $\bu = (u_1,\dots,u_n)$
  for $s$ such that $u_k(p) = p$ for all $i < k \leq j$ and $1 \leq p
  \leq s_{ij}$.  Since no lace diagram $\bw = (w_1,\dots,w_n)$ for $r$
  can satisfy these requirements, some $w_k$ is not dominated by $u_k$
  in the Bruhat order, which implies that $\phi_\bu(\Schub(\bw)) = 0$.
  We therefore get $\phi_s(Q_r) = \phi_\bu(Q_r) = 0$ which proves (i).
  
  For (ii), let $\bu$ be a fixed minimal lace diagram for $r$.  If
  $\bw$ is any minimal lace diagram for this orbit such that
  $\phi_\bu(\Schub(\bw)) \neq 0$, then since $w_i \leq u_i$ for all
  $i$ we must have $\bw = \bu$.  It therefore suffices to show that
  $\phi_\bu(\Schub(\bu)) = e(r)$, which follows from
  (\ref{eqn:schub_spec}) because $\phi_\bu$ maps each factor
  $\Schub_{u_i}(x^i;x^{i-1})$ to the product of the differences $(b_p
  - b_q)$ corresponding to strands of $\bu$ that cross between column
  $i-1$ and column $i$.  This finishes the proof.
\end{proof}


\section{Proof using Gr\"obner degeneration}
\label{sec:grobner}

In \cite{knutson.miller.ea:four} the closure of an orbit $r$ in the
space of quiver representations $V$ was degenerated into a union of
products of matrix Schubert varieties.  As a consequence of this, it
was proved \cite[Cor.~4.9]{knutson.miller.ea:four} that the Thom
polynomial $\tp_r$ can be written as a non-negative linear combination
of products of Schubert polynomials, indexed by minimal lace diagrams
for $r$.  In this section we give a new proof of the component formula
based on this fact.  The crucial observation is that a linear
combination of Schubert products can only be symmetric if all
coefficients are equal.

Consider any linear combination
\[ P = \sum_\bw c_\bw \, \Schub_{w_1}(x^1;x^0)\, \Schub_{w_2}(x^2;x^1)
   \cdots \Schub_{w_n}(x^n;x^{n-1})
\]
where the sum is over all minimal lace diagrams $\bw =
(w_1,\dots,w_n)$ for an orbit $r$.  Recall from the proof of
\reflemma{lemma:comp_sym} that if a divided difference operator
$\partial^i_j$ is evaluated on $P$, with $1 \leq i \leq n-1$ and $1
\leq j < e_i$, then the result is a linear combination of products
$\Schub(\bu) = \prod \Schub_{u_i}(x^i;x^{i-1})$ for lace diagrams $\bu
= (u_1,\dots,u_n)$, such that $u_i(j) < u_i(j+1)$ and $u_{i+1}^{-1}(j)
< u_{i+1}^{-1}(j+1)$.  Furthermore, the coefficient of $\Schub(\bu)$
is equal to $c_{\bu'} - c_{\bu''}$, where $\bu' = (u_1,\dots,u_i s_j,
u_{i+1},\dots,u_n)$ and $\bu'' = (u_1,\dots,u_i, s_j u_{i+1}, \dots,
u_n)$.  It follows that if $P$ is symmetric in all groups of variables
$x^i$, then for any minimal lace diagram $\bw = (w_1,\dots,w_n)$ such
that $w_i(j)>w_i(j+1)$ or $w_{i+1}^{-1}(j) > w_{i+1}^{-1}(j+1)$ we
have $c_\bw = c_{\bw'}$ where $\bw' = (w_1,\dots,w_i s_j, s_j
w_{i+1},\dots,w_n)$.  The transformation from $\bw$ to $\bw'$ is
illustrated by the following picture (of parts of the extended lace
diagrams):
\begin{equation} \label{eqn:trans}
\raisebox{-9pt}{\pic{50}{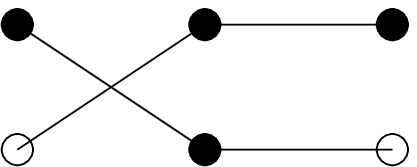}} 
\ \ \  \longleftrightarrow \ \ \ 
\raisebox{-9pt}{\pic{50}{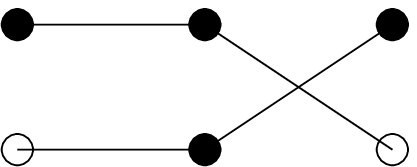}}
\end{equation}
Notice that this transformation can be applied to any lace diagram, as
long as the upper and middle dots are not in the extended part of the
diagram. 

We define the {\em left-most lace diagram\/} for the orbit $r$ as
follows.  Start with an empty diagram (with zero dots in each column).
Then for each $i = 0,1,\dots,n$, and each $j = n, n-1,\dots,i$ (in
this order) we add $r_{ij} - r_{i-1,j} - r_{i,j+1} + r_{i-1,j+1}$
strands starting at column $i$ and terminating at column $j$ to the
bottom of the diagram.  Notice that any left-most diagram is also
minimal.  The following picture shows an example of a left-most lace
diagram.
\[ \pic{50}{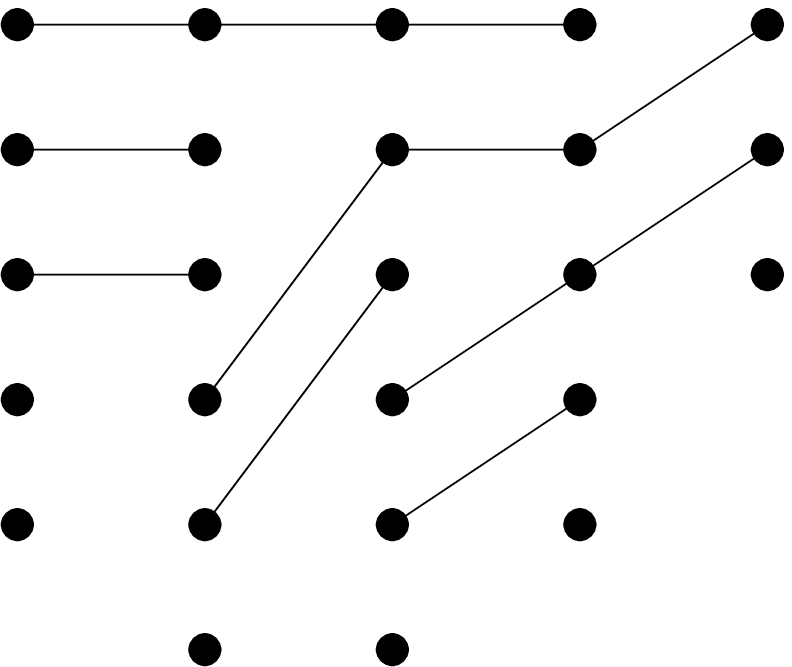} \]

\begin{prop} \label{prop:sym_crit}
  {\rm (i)} Two minimal lace diagrams are connected via the transformations
  (\ref{eqn:trans}) if and only if they represent the same orbit.
  
  {\rm (ii)} A linear combination $P = \sum c_\bw \Schub(\bw)$ over minimal
  lace diagrams $\bw$ for an orbit $r$ is symmetric in each group of
  variables $x^i$ if and only if all coefficients $c_\bw$ are equal.
\end{prop}
\begin{proof}
  It is enough to show that any minimal lace diagram can be converted
  to a left-most lace diagram using the transformations
  (\ref{eqn:trans}).  We give an explicit algorithm for doing this.
  
  Consider the strand which starts at the top dot of column $0$ in the
  lace diagram.  If this strand is not entirely in the top row of the
  diagram, we let $i$ be the first column where the strand contains a
  dot below the top row, and $k > 1$ the row number of this dot.  The
  line segment to this dot from the top dot of column $i-1$ must then
  cross the strand going through dot $k-1$ of column $i$, so these
  strands do not cross between column $i$ and column $i+1$.
  Furthermore, since these strands cannot terminate in the same
  column, the strand through dot $k-1$ of column $i$ continues to a
  dot of column $i+1$ which is not in the extended part of the
  diagram.  We can therefore use a transformation (\ref{eqn:trans}) to
  move the crossing one step to the right; in the resulting diagram,
  the strand starting at the top dot of the first column will now
  contain the $k-1$'st dot of column $i$.  By continuing to apply this
  method, we eventually reach a lace diagram in which the strand
  starting at the top dot of column $0$ is entirely in the top row.
  The same procedure is now carried out for the remaining strands that
  start at the first column, from top to bottom, then the strands
  starting at the second column, and so on; for each of these strands
  one ignores the previous strands which have already been moved to
  the correct position.  Finally, notice that since strands starting
  in the same column cannot cross each other, this algorithm will deal
  with the strands in the same order as they are added when a
  left-most lace diagram is constructed.  We conclude that the
  resulting lace diagram is left-most.
\end{proof}

\begin{proof}[Proof of \refthm{thm:ns_comp} using Gr\"obner
  degeneration]
  
  By \cite[Cor.~4.9]{knutson.miller.ea:four} and part (ii) of
  \refprop{prop:sym_crit}, the Thom polynomial $\tp_r$ is equal to a
  non-negative integer $c$ times the polynomial $Q_r$.  By the
  corollary to \reflemma{lemma:comp_sym} this says that
\[ \tp_r = c \sum_\bw \prod_{i=1}^n F_{w_i}(x^i;x^{i-1}) \]
where the sum is over all minimal lace diagrams for $r$.  Since each
Stanley symmetric function $F_{w_i}$ is an integral linear combination
of Schur polynomials, it follows that $c$ must divide all the quiver
coefficients $c_\lambda(r)$ for the orbit $r$.  To show that $c=1$ it
is therefore enough to find a quiver coefficient equal to one.

This can be done explicitly as follows.  For all $0 \leq i < j \leq n$
we let $R_{ij}$ be a rectangular partition with $r_{i+1,j} - r_{ij}$
rows and $r_{i,j-1} - r_{ij}$ columns, and we let $\lambda_i$ be the
Young diagram obtained by arranging the rectangles $R_{i-1,j}$ for $i
\leq j \leq n$ side by side from left to right.
\[ \lambda_i \ = \ \raisebox{-18pt}{\pic{70}{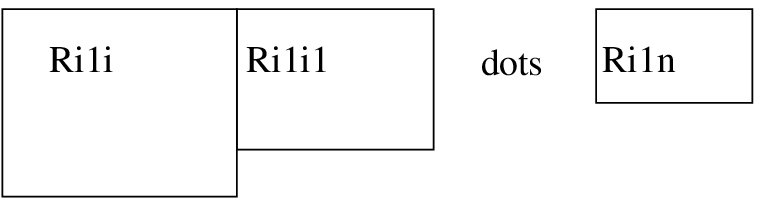}} \]
It then follows from the algorithm of \cite[\S 2.1]{buch.fulton:chern}
that $c_\lambda(r) = 1$ for the sequence of partitions $\lambda =
(\lambda_1,\dots,\lambda_n)$.
%
\end{proof}

\begin{remark}
  (a) With slightly more care, one can use the transformations
  (\ref{eqn:trans}) to prove that a linear combination $\sum c_\bw
  \Schub(\bw)$ over all lace diagrams $\bw$ for a given dimension
  vector is symmetric in each set of variables $x^i$ if and only if
  all coefficients corresponding to non-minimal lace diagrams are
  zero, and coefficients for minimal lace diagrams representing the
  same orbit are equal.
  
  (b) M.~Shimozono reports that the moves of (\ref{eqn:trans}) can
  also be used to prove that the components of the Gr\"obner
  degeneration of an orbit closure intersect in codimension two or
  higher.
\end{remark}


\section{Grothendieck classes of quiver varieties}
\label{sec:ktheory}

In \cite{buch:grothendieck} a formula for the Grothendieck class of a
quiver variety was proved, which generalizes (\ref{eqn:orig_qf}).
This formula can be interpreted as an expression for the structure
sheaf $\cO_{\wb r}$ of an orbit closure $\wb r$ in the
torus-equivariant Grothendieck ring of the representation space $V$
\cite{knutson.miller.ea:four,fulton:notes}.  It has the form
\begin{equation} \label{eqn:orig_kqf}
  [\cO_{\wb r}] = \sum_\lambda c_\lambda(r)\, G_{\lambda_1}(x^1;x^0) \,
  G_{\lambda_2}(x^2;x^1) \cdots G_{\lambda_n}(x^n;x^{n-1})
\end{equation}
where $G_{\lambda_i}$ denotes the double stable (Laurent) Grothendieck
polynomial for the partition $\lambda_i$ (see \cite[\S
2]{buch:alternating} for notation).  The sequences $\lambda$ of
partitions in this formula all satisfy that the sum $\sum |\lambda_i|$
of the weights is greater than or equal to the expected codimension
$d(r)$.  The cohomological quiver coefficients of (\ref{eqn:orig_qf})
are the subset of the coefficients $c_\lambda(r)$ in
(\ref{eqn:orig_kqf}) for which $\sum |\lambda_i| = d(r)$.  It was
conjectured in \cite{buch:grothendieck} that the $K$-theoretic quiver
coefficients have signs which alternate with codimension, that is
$(-1)^{\sum |\lambda_i| - d(r)}\, c_\lambda(r) \geq 0$.

This conjecture was proved in \cite{buch:alternating} by giving a
$K$-theoretic generalization of the component formula.  E.~Miller has
found a different proof of this formula \cite{miller:alternating}.
The $K$-theoretic component formula has the form
\begin{equation} \label{eqn:k_comp}
  [\cO_{\wb r}] = \sum_\bw\, (-1)^{\sum \ell(w_i)-d(r)} \,
  \Groth_{w_1}(x^1;x^0)\, \Groth_{w_2}(x^2;x^1) \cdots
  \Groth_{w_n}(x^n;x^{n-1})
\end{equation}
where $\Groth_{w_i}$ is the (Laurent) Grothendieck polynomial of
Lascoux and Sch\"utzen\-berger \cite{lascoux.schutzenberger:structure,
  lascoux:anneau}, and the sum is over certain lace diagrams $\bw =
(w_1,\dots,w_n)$ called {\em KMS-factorizations\/} for the orbit $r$.
These lace diagrams can be defined as certain factorizations of the
Zelevinsky permutation of \cite{knutson.miller.ea:four}.  In this
final section we explain how the methods of the present paper can be
used to give a concrete description of the KMS-factorizations
associated to a given orbit $r$.

Let $P = \sum c_\bw \Groth(\bw)$ be a linear combination of products
of Grothendieck polynomials $\Groth(\bw) = \prod
\Groth_{w_i}(x^i;x^{i-1})$ for all lace diagrams $\bw$ for the
dimension vector $(e_0,\dots,e_n)$.  The arguments of section
\ref{sec:grobner} can be generalized to show that $P$ is symmetric in
each set of variables $x^i$ if and only if the following conditions
are satisfied:

(I) The coefficient $c_\bw$ of a lace diagram $\bw = (w_1,\dots,w_n)$
is non-zero only if $w_1^{-1}(j) < w_1^{-1}(j+1)$ for all $1 \leq j <
e_0$ and $w_n(j) < w_n(j+1)$ for all $1 \leq j < e_n$.

(II) For every lace diagram $\bu = (u_1,\dots,u_n)$ and integers $1
\leq i < n$ and $1 \leq j < e_i$ such that $u_i(j) < u_i(j+1)$ and
$u_{i+1}^{-1}(j) < u_{i+1}^{-1}(j+1)$, we have $c_{\bu'} = c_{\bu''} =
-c_{\bu'''}$, where $\bu' = (u_1,\dots,u_i s_j, u_{i+1},\dots,u_n)$, $\bu'' =
(u_1,\dots,u_i,s_j u_{i+1},\dots,u_n)$, and $\bu''' = (u_1,\dots,u_i
s_j, s_j u_{i+1}, \dots, u_n)$.  

It follows easily from the definition of KMS-factorizations given in
\cite{buch:alternating} that any KMS-factorization $\bw$ satisfies the
requirement of (I), and that each of the lace diagrams $\bu'$,
$\bu''$, and $\bu'''$ of (II) are KMS-factorizations for $r$ if and
only if all three are KMS-factorizations for $r$ (see the remark at
the end of \cite[\S 6]{buch:alternating}).  Since this is sufficient
to prove the description of KMS-factorizations presented here, we will
skip the proof of the above classification of symmetric linear
combinations of products of Grothendieck polynomials.  Notice that the
transformation on lace diagrams corresponding to (II) can be pictured
as follows.
\begin{equation} \label{eqn:k_trans}
\raisebox{-9pt}{\pic{50}{tx_.eps}} 
\ \ \  \longleftrightarrow \ \ \ 
\raisebox{-9pt}{\pic{50}{t_x.eps}}
\ \ \  \longleftrightarrow \ \ \ 
\raisebox{-9pt}{\pic{50}{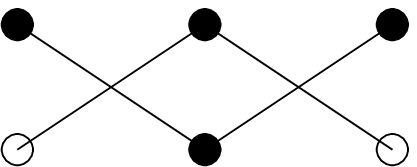}}
\end{equation}

We also need to know that a minimal lace diagram for an orbit $r$ is
a KMS-factorization for this orbit and no other orbits.  Again, this
fact is immediate from the definition of KMS-factorizations
\cite{buch:alternating}.

\begin{thm}
  The KMS-factorizations for an orbit $r$ are exactly the lace
  diagrams that can be obtained by applying a series of
  transformations (\ref{eqn:k_trans}) to the left-most lace diagram
  for $r$.
\end{thm}
\begin{proof}
  Let $\bw = (w_1,\dots,w_n)$ be any KMS-factorization for $r$.  By
  applying a series of transformations (\ref{eqn:k_trans}) to $\bw$,
  each replacing an occurrence of the first or the third diagram of
  (\ref{eqn:k_trans}) with the middle diagram, one arrives at a
  KMS-factorization $\bw' = (w_1',\dots,w_n')$ in which only the
  middle situation of (\ref{eqn:k_trans}) can be found.  It is enough
  to prove that $\bw'$ is a minimal lace diagram.  In fact, if this is
  true then $\bw'$ must be the left-most diagram for $r$, since the
  algorithm in the proof of \refprop{prop:sym_crit} will not change
  this diagram.
  
  If $\bw'$ has a crossing outside the extended part of the diagram,
  say between column $i-1$ and column $i$, then one can find $1 \leq j
  < e_i$ such that $w_i(j) > w_i(j+1)$.  By (I) this implies that $i <
  n$.  Since the first and third situations of (\ref{eqn:k_trans})
  cannot occur, the two crossing strands must both terminate at column
  $i$, which implies that $\bw'' = (w'_1, \dots, w'_i s_j, s_j
  w'_{i+1}, \dots, w'_n)$ is not a lace diagram.  On the other hand,
  (II) requires $\bw''$ to be a KMS-factorization, a contradiction.
  
  We conclude from this that every crossing of $\bw'$ must involve a
  line segment in the extended part of the diagram, which extends the
  right end of a strand.  In particular, two strands can cross at most
  once, and not at all if they terminate at the same column.
  
  It remains to show that no two crossing strands of $\bw'$ can start
  at the same column.  Assume for contradiction that a strand starting
  at dot $j$ of column $i$ crosses another starting at dot $k$ of
  column $i$, where $j < k$.  Assume also that $k-j$ is minimal with
  these properties.  Since all crossings involve line segments
  extending the right end of a strand, it follows that the strand
  starting at dot $j$ is shorter than the strand starting at dot $k$.
  Furthermore, if $j+1 < k$ then the strand containing dot $j+1$ of
  column $i$ must start at this dot; otherwise it would cross the left
  side extension of the strand starting at dot $j$.  Since the strand
  starting at dot $j+1$ is either longer than the strand starting at
  dot $j$ or shorter than the strand starting at dot $k$, the
  minimality of $k-j$ forces $k = j+1$.  Now a series of the moves
  (\ref{eqn:trans}), from right to left, will move the crossing of the
  two strands so that it occurs between columns $i$ and $i+1$.  But
  this is again impossible: (I) implies that $i > 0$, after which (II)
  can be used to produce a KMS-factorization which is not a lace
  diagram.  This contradiction shows that $\bw'$ is a minimal lace
  diagram, which concludes the proof.
\end{proof}


\bibliographystyle{amsplain}
\bibliography{warsawa}

\end{document}